\documentclass[11pt,a4paper]{amsart}
\usepackage[centering,left=3cm,right=3cm]{geometry}
\usepackage{amsfonts,amscd,amssymb,amsmath,amsthm,mathrsfs,xcolor,lscape,longtable}
\usepackage{microtype}
\usepackage[Algorithm]{algorithm}
\usepackage[noend]{algpseudocode}
\usepackage[colorlinks,linkcolor=blue,anchorcolor=blue,citecolor=blue,backref=page]{hyperref}
\usepackage{url}
\usepackage{enumitem}
\usepackage{color}
\usepackage{graphics,epsfig}
\usepackage{float}
\usepackage{epstopdf}
\usepackage[utf8]{inputenc}
\usepackage{hyperref}
\usepackage{calc, amscd}
\usepackage[T1]{fontenc}
\usepackage{pdfpages}
\usepackage{mathtools} 
\usepackage{scalerel,stackengine}
\hypersetup{breaklinks=true}
\usepackage{multirow,pbox,lipsum}
\usepackage[thinlines]{easytable}
\usepackage{todonotes}
\usepackage{phaistos}

 \usepackage[all,cmtip]{xy}\usepackage{tikz}\usetikzlibrary{arrows}
\usetikzlibrary{matrix}\usetikzlibrary{calc}
\usepackage{tikz-cd,xr,multicol,microtype}

\newtheorem{thm}{Theorem}

\newtheorem{defn}[thm]{Definition}
\newtheorem{cor}[thm]{Corollary}

\theoremstyle{definition}
\newtheorem{rmk}{Remark}
\theoremstyle{definition}

\numberwithin{equation}{section}

\newcommand{\bb}[1]{\mathbb{#1}}

\newcommand{\C}{\bb{C}}
\newcommand{\Q}{\bb{Q}}
\newcommand{\R}{\bb{R}}
\newcommand{\Z}{\bb{Z}}

\newcommand{\rU}{\mathrm{U}}
\newcommand{\rB}{\mathrm{B}}

\newcommand{\G}{\Gamma}
\newcommand{\Sp}{\mathrm{Sp}}

\newcommand{\p}{\mathbb{P}}

\newcommand{\GL}{\mathrm{GL}}

\begin{document}
\date{\today}
\title[$\Sp(6)$  hypergeometric groups]{commentary on $\Sp(6)$ hypergeometric groups}
\author{Jitendra Bajpai}
\address{Max Planck Institute for Mathematics, Vivatsgasse 7, 53111 Bonn, Germany}
\email{jitendra@mpim-bonn.mpg.de}
\subjclass[2010]{Primary: 22E40;  Secondary: 32S40;  33C80}  
\keywords{Hypergeometric groups, monodromy representations, symplectic groups}


\begin{abstract}
We study the examples mentioned in~\cite[Tables~A\&~C]{BDSS} and establish the arithmeticity of four examples of symplectic hypergeometric groups of degree six. 

Following~\cite{BDSS} we know that there are 458 inequivalent symplectic hypergeometric groups of degree six, and combining the results of this article with the work of~\cite{BDN, BDSS, SV}, we now know that  at least 384 are arithmetic and at least 63 are thin whereas the arithmeticity and thinness of remaining 11 examples are
still unknown.
\end{abstract}

\maketitle


\section{Introduction}\label{sec:intro}
A hypergeometric differential equation of order $n$ is an ordinary differential equation of order $n$ with three regular singular points. It is defined on the thrice punctured Riemann sphere $\p^{1}(\C)\backslash \{ 0,1,\infty\}$. Let $\theta = z \frac{d}{dz}$ and the parameters $$\alpha = (\alpha_1, \ldots , \alpha_n) , \beta = ( \beta_1, \ldots , \beta_n)  \in \C^{n}.$$ We define the \textit{hypergeometric differential equation} of order $n$ by 
\begin{equation}\label{hde}
[z(\theta + \alpha_1) \cdots (\theta + \alpha_n) - (\theta+\beta_1 -1)\cdots (\theta+\beta_{n} -1 )] u(z) =0\,.
\end{equation}
This has $n$ linearly independent solutions which can be explicitly expressed as hypergeometric functions of type ${}_{n} F_{n-1}$ around any point $z \in \p^{1}(\C)\backslash \{ 0,1,\infty\}$. For $\alpha $ and $\beta$, we define 

$${}_{n} F_{n-1}(\alpha_1,\ldots, \alpha_n; \beta_1, \ldots, \beta_{n-1} | z)= \sum_{k=0}^{\infty} \frac{(\alpha_1)_{k} \ldots (\alpha_n)_{k} }{(\beta_1)_{k} \ldots (\beta_{n-1})_{k}} \frac{z^{k}}{k!}\,,$$ where $(\alpha)_{k}= \frac{\Gamma(\alpha+k)}{\Gamma(\alpha)}$. Then a system of $n$ linearly independent solutions of~\eqref{hde} is defined by the hypergeometric functions of type ${}_{n} F_{n-1}$.

The fundamental group $\pi_1$ of $\p^{1}(\C)\backslash \{ 0,1,\infty\}$ acts on the (local) solution space of the hypergeometric equation \eqref{hde}, and we get the monodromy representation $\rho:\pi_1\longrightarrow \GL(V)$ where $V$ is the $n$-dimensional solution space of the differential equation~\eqref{hde} on a small neighbourhood of a point in $\p^{1}(\C)\backslash \{ 0,1,\infty\}$. 

The subgroup $\rho(\pi_1)$ of $\GL(V)$ is said to be the \textit{monodromy group} of the hypergeometric differential equation~\eqref{hde}, which is now called as the \textit{hypergeometric group} associated to the parameters $\alpha , \beta \in \C^{n}.$ 

Levelt (cf.~\cite[Theorem 3.5]{BH}) showed that, if $\alpha_j-\beta_k\notin\Z$ for all $1\le j,k\le n$, then there exists a basis of the solution space of the hypergeometric equation with respect to which the hypergeometric group corresponding to the parameters $\alpha , \beta  \in \C^{n}$ is the subgroup of $\GL_n(\C)$ generated by the companion matrices $A$ and $B$ of the polynomials
\[f(x)=\prod_{j=1}^n(x-e^{2\pi i\alpha_j}),\quad g(x)=\prod_{j=1}^n(x-e^{2\pi i\beta_j})\] respectively, and the monodromy is defined by  \begin{equation}\label{eq:monodromy}
g_\infty\mapsto A , \quad  g_0\mapsto B^{-1} , \quad  g_1\mapsto A^{-1}B , 
\end{equation}
where $g_0, g_1, g_\infty$ are, respectively, the loops around $0,1,\infty$, which generate the fundamental group of $\mathbb{P}^1(\C)\backslash\{0,1,\infty\}$ modulo the relation $g_\infty g_1 g_0=1$. Any other hypergeometric group having the same parameters is a conjugate of this one. Note that the condition $\alpha_j-\beta_k\notin\Z$ for all $1\le j,k\le n$ ensures that the polynomials $f$ and $g$ do not have any common root. 

Let $\Gamma(f,g)$ denote the hypergeometric group generated by the companion matrices of the polynomials $f,g$; as said above, it is a subgroup of $\GL_n(\C)$. We consider the cases where the coefficients of $f,g$ are integers with $f(0)=\pm1$, $g(0)=\pm1$ (for example, take $f,g$ as products of cyclotomic polynomials); in these cases, $\Gamma(f,g)\subset\GL_n(\Z)$. In addition, we assume that $f,g$ form a primitive pair~\cite[Definition 5.1]{BH}, are self-reciprocal and do not have any common root.
 
Beukers and Heckman \cite[Theorem 6.5]{BH} have determined the Zariski closures $\mathrm{G}$ of the hypergeometric groups $\Gamma(f,g)$, which can be described as follows.

\begin{itemize}
\item If $n$ is even and $f(0)=g(0)=1$, then the hypergeometric group $\Gamma(f,g)$ preserves a non-degenerate integral symplectic form $\Omega$ on $\Z^n$ and $\Gamma(f,g)\subset\Sp_\Omega(\Z)$ is Zariski dense, that is, $\mathrm{G}=\Sp_\Omega$.

\item If $\Gamma(f,g)$ is {\it infinite} and $\frac{f(0)}{g(0)}=-1$, then $\Gamma(f,g)$ preserves a non-degenerate integral quadratic form $\mathrm{Q}$ on $\Z^n$ and $\Gamma(f,g)\subset\mathrm{O}_\mathrm{Q}(\Z)$ is Zariski dense, that is, $\mathrm{G}=\mathrm{O}_\mathrm{Q}$.

\item It follows from~\cite[Corollary 4.7]{BH} that $\Gamma(f,g)$ is {\it finite} if and only if either $\alpha_1<\beta_1<\alpha_2<\beta_2<\cdots<\alpha_n<\beta_n$ or $\beta_1<\alpha_1<\beta_2<\alpha_2<\cdots<\beta_n<\alpha_n$. In this case, we say that the roots of $f$ and $g$ {\it interlace} on the unit circle.
\end{itemize}
 \begin{defn}\label{defn:at} 
 We call a hypergeometric group $\Gamma(f,g) \subseteq \mathrm{G}(\Z)$ {\it arithmetic} if it is of {\it finite} index in $\mathrm{G}(\Z)$, and {\it thin} if it has {\it infinite} index in $\mathrm{G}(\Z)$, where $\mathrm{G}$ is the Zariski closure of $\Gamma(f,g)$ inside $\GL_n(\C)$.
\end{defn}

We take this opportunity to remind the readers that this is not the most general definition of an arithmetic  group. However, the groups $\Gamma(f,g)$ under consideration are simply subgroups of $\GL_n(\Z)$ and therefore the above definition to define an arithmetic group is quite natural. For quick introduction on the theory of arithmetic groups we refer the interested reader  to~\cite{Borel}.

Sarnak's question~\cite{Sa14} about classifying the pairs of polynomials $f,g$ for which the associated hypergeometric group $\Gamma(f,g)$ is arithmetic or thin, has witnessed many interesting developments over the past 10 years, and this article,  with no exception, is a small addition to these developments.  

More precisely, this article is an extension of the work carried out by the author and his collaborators in~\cite{BDSS} and~\cite{BDN} about determining the arithmeticity and thinness of symplectic hypergeometric groups of degree six. 

In~\cite{BDSS}, the authors showed that there are in total 458 pairs of polynomials $f, g$ (up to scalar shifts) which define symplectic hypergeometric groups $\Gamma(f,g)$ of degree six. For 211 of them, the absolute value of the leading coefficients of the difference polynomials $f-g$ are at most 2 and therefore the arithmeticity of these 211 examples simply follows from Theorem~1.1 of~\cite{SV}. The arithmeticity of 164 examples follows from Theorems~2 and~3 of~\cite{BDSS}. Note that among these 164 arithmetic examples, arithmeticity of one group, that is the example~\cite[Table B-145]{BDSS}, was proved in~\cite{DFH} by computing explicitly the index of this particular group inside $\Sp_6(\Z)$. In~\cite{BDSS}, the arithmeticity of these 164 examples were shown by  finding an element $\gamma$ satisfying the hypotheses of Proposition~1 from~\cite[Page~260]{BDSS}. Further, in~\cite{BDN} the authors manage to show the arithmeticity of 5 more examples by finding an element $\gamma$ satisfying the hypotheses of this Proposition~1 in~\cite{BDSS}. Moreover, the authors also proved the thinness of 63 examples in~\cite{BDN} by using a version of the well-known ``ping-pong lemma'' from geometric group theory. At this stage, there were in total 15 examples left whose arithmeticity and thinness were still unknown. 

For these 15 remaining cases, in~\cite{BDN}, the authors also pointed out that by all means our approach to play ping-pong will remain inconclusive. Moreover, there are 6 examples where the arithmeticity can not be concluded by Prop.~1 of~\cite{BDSS}, since it will be impossible to find an element $\gamma$ which satisfies the hypotheses of Prop.~1. This is due to the gcd of the coordinates of vector $v=(A^{-1}B -I)e_6$, where $e_6=(0,0,0,0,0,1)$, which is larger than 2 in these 6 cases, see the entries of the last column in Tables~\ref{ant} and~\ref{open}. These 6 cases are labelled as A-1, C-1, C-10, C-42, C-59, C-61 in the Tables~\ref{ant} and~\ref{open} below, to which we focused on in this article to try to prove their arithmeticity. At the end, we are successful in establishing the arithmeticity of 4 out of these 6 examples by using the method used to prove one of the main results of~\cite[Thm.~1.2]{SV}. More precisely, the work of this article provide the following result. 
\begin{thm}\label{thm:main}
The hypergeometric groups associated to the four pairs of parameters $\alpha, \beta$ of Table~\ref{ant} are arithmetic.
\end{thm}

{\begin{center}
\begin{table}[htbp]
\caption{Arithmetic hypergeometric groups in $\Sp(6)$}\label{ant}
\scriptsize\renewcommand{\arraystretch}{2}
\begin{tabular}{|c|c|c|c|c|c|}
\hline
 label &$\alpha$ & $\beta$  & $v$     \\                                                                                                                                                        
\hline
\hline
 C-1 &$\big(0,0,0,0,\frac{1}{2},\frac{1}{2}\big)$ & $\big(\frac{1}{3},\frac{1}{3},\frac{2}{3},\frac{2}{3},\frac{1}{6},\frac{5}{6}\big)$  & $\left(-3\,,-3\,,3\,,-3\,,-3\,,0\right)$  \\
\hline
 C-10 &$\big(0,0,0,0,\frac{1}{3},\frac{2}{3}\big)$& $\big(\frac{1}{9},\frac{2}{9},\frac{4}{9},\frac{5}{9},\frac{7}{9},\frac{8}{9}\big)$  &$\left(-3\,,3\,,-3\,,3\,,-3\,,0\right)$ \\
\hline
C-42 &$(0,0,\frac{1}{4},\frac{1}{4},\frac{3}{4},\frac{3}{4})$& $(\frac{1}{3},\frac{2}{3},\frac{1}{12},\frac{5}{12},\frac{7}{12},\frac{11}{12})$  &$\left(-3\,,3\,,-3\,,3\,,-3\,,0\right)$ \\
\hline
 C-59 &$(0,0,\frac{1}{12},\frac{5}{12},\frac{7}{12},\frac{11}{12})$& $\big(\frac{1}{3},\frac{2}{3},\frac{1}{4},\frac{3}{4},\frac{1}{4},\frac{3}{4}\big)$ & $\left(-3\,,-3\,,0\,,-3\,,-3\,,0\right)$  \\
\hline
\end{tabular}
\end{table}
\end{center}
}

\begin{rmk}
We will label the hypergeometric groups discussed in this article according to how they appear in Table A and Table C of~\cite{BDSS}: ``X-Y'' will represent the entry Y from Table X of~\cite{BDSS}.
\end{rmk}

Now, following~\cite[Thm.~1.1]{SV}, \cite[Thms.~2 \&~3]{BDSS}, \cite[Thms.~2,\,3,\,5 \&~6]{BDN}, along with Theorem~\ref{thm:main} of this article, we describe the current state of the problem about determining arithmeticity and thinness of symplectic hypergeometric groups of degree six in the form of following corollary.
\begin{cor}\label{co:total}
There are $458$ degree six symplectic hypergeometric groups up to equivalence. Among them, at least $384$ are arithmetic and $63$ are thin.
\end{cor}

At the end, we are only left with 11 examples of degree six hypergeometric groups whose arithmeticity and thinness are still \emph{unknown}. We list these remaining 11 cases in Table~\ref{open} below.

{\renewcommand{\arraystretch}{2}
\begin{table}[h]
\small
\caption{Open cases}\label{open}
\addtolength{\tabcolsep}{-1pt}
\newcounter{rownum-0}
\setcounter{rownum-0}{0}
\centering
\begin{tabular}{|c|c|c|c|c|c|}
\hline
 Label & $\alpha$ & $\beta$  & $v$ \\                                                            
\hline
\hline
A-15 & $\left(0,0,0,0,0,0\right)$& $\big(\frac{1}{3},\frac{1}{3},\frac{1}{3},\frac{2}{3},\frac{2}{3},\frac{2}{3}\big)$   &  $\left(-9\,,9\,,-27\,,9\,,-9\,,0\right)$ \\  
\hline
A-16 & $\left(0,0,0,0,0,0\right)$& $\big(\frac{1}{3},\frac{1}{3},\frac{2}{3},\frac{2}{3},\frac{1}{4},\frac{3}{4}\big)$  &  $\left(-8\,,11\,,-24\,,11\,,-8\,,0\right)$ \\ 
\hline
A-21 &$\left(0,0,0,0,0,0\right)$ & $\big(\frac{1}{3},\frac{2}{3},\frac{1}{5},\frac{2}{5},\frac{3}{5},\frac{4}{5}\big)$   & $\left(-8\,,12\,,-23\,,12\,,-8\,,0\right)$ \\  
\hline
C-9 & $\left(0,0,0,0,\frac{1}{3},\frac{2}{3}\right)$&$\left(\frac{1}{7},\frac{2}{7},\frac{3}{7},\frac{4}{7},\frac{5}{7},\frac{6}{7}\right)$ & $\left(-4\,,2\,,-3\,,2\,,-4\,,0\right)$\\
 \hline
C-31 & $\left(0,0,0,0,\frac{1}{6},\frac{5}{6}\right)$&$\left(\frac{1}{3},\frac{2}{3},\frac{1}{5},\frac{2}{5},\frac{3}{5},\frac{4}{5}\right)$ & $\left(-7\,,8\,,-17\,,8\,,-7\,,0\right)$\\
 \hline
C-32 & $\left(0,0,0,0,\frac{1}{6},\frac{5}{6}\right)$&$\left(\frac{1}{4},\frac{3}{4},\frac{1}{12},\frac{5}{12},\frac{7}{12},\frac{11}{12}\right)$& $\left(-5\,,11\,,-14\,,11\,,-5\,,0 \right)$\\
\hline
C-47 &$\left(0,0,\frac{1}{5},\frac{2}{5},\frac{3}{5},\frac{4}{5}\right)$ & $\left(\frac{1}{2},\frac{1}{2},\frac{1}{3},\frac{1}{3},\frac{2}{3},\frac{2}{3}\right)$ & $\left(-5\,,-8\,,-10\,,-8\,,-5\,,0\right)$\\
\hline
C-51 & $\left(0,0,\frac{1}{6},\frac{1}{6},\frac{5}{6},\frac{5}{6}\right)$ & $\left(\frac{1}{2},\frac{1}{2},\frac{1}{12},\frac{5}{12},\frac{7}{12},\frac{11}{12}\right)$& $\left(-6\,,8\,,-8\,,8\,,-6\,,0 \right)$ \\
\hline
C-55 &$\left(0,0,\frac{1}{8},\frac{3}{8},\frac{5}{8},\frac{7}{8}\right)$ & $\left(\frac{1}{2},\frac{1}{2},\frac{1}{12},\frac{5}{12},\frac{7}{12},\frac{11}{12}\right)$  & $\left( -4\,,1\,,2\,,1\,,-4\,,0 \right)$\\
\hline
C-60 & $\left(\frac{1}{3},\frac{1}{3},\frac{1}{3},\frac{2}{3},\frac{2}{3},\frac{2}{3}\right)$&$\left(\frac{1}{6},\frac{1}{6},\frac{1}{6},\frac{5}{6},\frac{5}{6},\frac{5}{6}\right)$ & $\left(6\,,0\,,14\,,0\,,6\,,0 \right)$\\
\hline
C-61 & $\left(\frac{1}{3},\frac{1}{3},\frac{1}{3},\frac{2}{3},\frac{2}{3},\frac{2}{3}\right)$ & $\left(\frac{1}{9},\frac{2}{9},\frac{4}{9},\frac{5}{9},\frac{7}{9},\frac{8}{9}\right)$ & $\left(3\,,6\,,6\,,6\,,3\,,0\right)$ \\
\hline
\end{tabular}
\end{table}
}
Note that 9 out of these 11 open cases, are still available to try to show the arithmeticity by finding such an element $\gamma$, if we can. However, the method which we adapt in this article to show the arithmeticity of 4 examples, can simply be attempted to establish the arithmeticity of all the remaining 11 cases.

We remind the reader that the thinness of the 63 examples mentioned in Corollary~\ref{co:total} above was achieved by adapting the approach of Brav and Thomas~\cite{BT} from dimension $n=4$ to $n=6$. Note that the first examples of higher rank thin hypergeometric groups were found in~\cite{BT} in the case of symplectic hypergeometric groups $\Gamma(f, g) \subset \Sp_4(\Z)$ by playing ping-pong.

\section{Preliminaries}\label{se:basic}

Let $f, g$ be a pair of degree $6$ polynomials that are products of cyclotomic polynomials, do not have any common roots, form a primitive pair (that is, there do not exist polynomials $f_1,g_1\in\Z[x]$ so that $f(x)=f_1(x^k), g(x)=g_1(x^k)$ for $k\ge 2$), and $f(0)=g(0)=1$. Observe that $f,g$ are self-reciprocal and the product of all the roots of $f$, as well as of $g$, is $1$. These conditions ensure that the corresponding hypergeometric group $\Gamma(f,g)$ preserves a nondegenerate symplectic form $\Omega$, and $\Gamma(f,g)$ is a Zariski dense subgroup of the corresponding symplectic group $\Sp_\Omega$ (cf.~\cite[Theorem 6.5]{BH}). Since the polynomials $f,g$ have integer coefficients and their constant terms are $1$, it follows that  $\Gamma(f,g)\subseteq\Sp_\Omega(\Z)$. 

For given polynomials $$f(x)=x^6+a_{5}x^{5}+\cdots+a_1x+1, \quad  g(x)=x^6+b_{5}x^{5}+\cdots+b_1x+1,$$ the corresponding companion matrices are

\resizebox{\textwidth}{!}{
\parbox{\textwidth}{
\begin{align*}
A=\left(\begin{array}{cccccccc} 
0&0&0&0&0&-1\\
1&0&0&0&0&-a_1\\
0&1&0&0&0&-a_{2}\\
0&0&1&0&0&-a_{3}\\
0&0&0&1&0&-a_{3}\\
0&0&0&0&1&-a_{5}
\end{array} \right),\quad
B=\left(\begin{array}{ccccccc} 
0&0&0&0&0&-1\\
1&0&0&0&0&-b_1\\
0&1&0&0&0&-b_{2}\\
0&0&1&0&0&-b_{3}\\
0&0&0&1&0&-b_{3}\\
0&0&0&0&1&-b_{5}
\end {array} \right)
\end{align*}}}
respectively, and  
$$C=A^{-1}B=\left(\begin{array}{cccccccc}
1&0&0&0&0&a_1-b_1\\ 
0&1&0&0&0& a_2-b_2\\ 
0&0&1&0&0&a_3-b_3\\ 
0&0&0&1&0&a_4-b_4\\
0&0&0&0&1&a_5-b_5\\ 
0&0&0&0&0&1\end {array}\right).$$

\subsection{Structure of the unipotent groups}
We will now briefly describe the structure of the unipotent groups corresponding to the roots of $\Sp_6(\Omega)$. Let $\mathcal{D}=\{ \epsilon_1, \epsilon_2, \epsilon_3, \epsilon_3^{*}, \epsilon_2^{*},\epsilon_1^{*}\}$ be a basis of $\Q^{6}$ over $\Q$, such that the matrix form $\Omega'$ of $\Omega$ with respect to this basis, is 

\resizebox{\textwidth}{!}{
\parbox{\textwidth}{
\begin{align*}
\Omega=\left(\begin{array}{cccccc} 
0&0&0&0&0&\lambda_1\\ 
0&0&0&0&\lambda_2&0\\ 
0&0&0&\lambda_3&0&0\\ 
0&0&-\lambda_3&0&0&0\\ 
0&-\lambda_2&0&0&0&0\\ 
-\lambda_1&0&0&0&0&0\end {array} \right)
\end{align*}
}}
with $\Omega(\epsilon_i, \epsilon_i^{*})=\lambda_i \in \Q^{*}$ for $1 \leq i \leq 3$. Let $T$ be the maximal torus in $\Sp_6(\Omega)$ given by the group of diagonal matrices, 

{\resizebox{\textwidth}{!}{
\parbox{\textwidth}{
\begin{align*}
T=\left\{ \left(\begin{array}{cccccc} 
t_1&0&0&0&0&0\\ 
0&t_2&0&0&0&0\\ 
0&0&t_3&0&0&0\\ 
0&0&0&t_3^{-1}&0&0\\ 
0&0&0&0&t_2^{-1}&0\\ 
0&0&0&0&0&t_1^{-1}\end {array} \right) | t_i \in \Q^{*}, \text{for} 1\leq i \leq 3\right\}
\end{align*}}}}
then $T$ defines a root system and, as $T$ is a $\Q$-split torus, this root system is also a $\mathbb{Q}$-root system for $\Sp_6$. More precisely, we may now define the root system $\Phi:=\Phi(T)$ for $\Sp_6(\Omega)$ by simply taking, for $1\leq i \leq 3,$ $t_i$ be the character of $T$ defined by $$\left(\begin{array}{cccccc} 
t_1&0&0&0&0&0\\ 
0&t_2&0&0&0&0\\ 
0&0&t_3&0&0&0\\ 
0&0&0&t_3^{-1}&0&0\\ 
0&0&0&0&t_2^{-1}&0\\ 
0&0&0&0&0&t_1^{-1}\end {array} \right) \mapsto t_i  \,.$$ Then the set of all roots is defined by $\Phi=\Phi^{+} \cup \Phi^{-}$, the union of the set of all positive and negative roots denoted by $\phi^{+}$ respectively $\Phi^{-}$. Hence, for a fix set of simple roots $\Pi=\{ t_1 t_2^{-1}, t_2 t_3^{-1}, t_3^{2}\}$, we have
\begin{align*}
&\Phi^{+} =\left\{t_1^2, t_1 t_2, t_1 t_3, t_1 t_3^{-1}, t_1 t_2^{-1}, t_2^2, t_2 t_3, t_2 t_3^{-1},t_3^2 \right\}\,, \\
&\Phi^{-} =\left\{t_1^{-2}, t_1^{-1} t_2^{-1}, t_1^{-1} t_3^{-1}, t_1^{-1} t_3, t_1^{-1} t_2, t_2^{-2}, t_2^{-1} t_3^{-1}, t_2^{-1} t_3,t_3^{-2} \right\}\,, 
\end{align*}
and therefore $t_1^{2}, t_1 t_2$ are the highest respectively the second highest roots in $\Phi^{+}$.
 We can now describe that the unipotent groups $U_{t_1^{2}}$ and $U_{t_1 t_2}$ corresponding to the highest and second highest roots as follows:
 
\resizebox{\textwidth}{!}{
\parbox{\textwidth}{
 \begin{align*}
& U_{t_1^{2}} =\left\{ \left(\begin{array}{rrrrrrr} 
1&0&0&0&0&{\it y}\\ 
0&1&0&0&0&0\\ 
0&0&1&0&0&0\\ 
0&0&0&1&0&0\\ 
0&0&0&0&1&0\\ 
0&0&0&0&0&1\end {array} \right)| y\in \Q\right\},\quad
U_{t_1 t_2} = \left\{\left(\begin{array}{rrrrrrrr} 
1&0&0&0&{\it x}&0\\ 
0&1&0&0&0&{\frac{\lambda_1}{\lambda_2}x}\\ 
0&0&1&0&0&0\\ 
0&0&0&1&0&0\\ 
0&0&0&0&1&0\\ 
0&0&0&0&0&1\end {array} \right)| x\in \Q, \,\,\lambda_1, \lambda_2 \, \text{as in} \, M\right\}.
\end{align*}}}

\subsection{Methodology}\label{se:meth} Note that to prove the arithmeticity of four examples in Table~\ref{ant}, we follow the method used to prove Theorem~1.2 of~\cite{SV}. However, the computations  are quite intricate and involved.

We start the process by computing the symplectic form $\Omega'$ (up to scalar multiples) preserved by the $\Gamma:=\Gamma(f,g)=\langle A, B \rangle$. Then, we show that there exists a basis  $\mathcal{D}=\{\epsilon_1,\epsilon_2, \epsilon_3, \epsilon_3^{*}, \epsilon_2^{*},\epsilon_1^{*} \}$ of $\Q^{6}$  such that the matrix form of $\Omega'$, with respect to this new basis $\mathcal{D}$, denoted by $\Omega''=X^{t}\Omega' X$ is anti-diagonal. Here, $X$ denote the change of basis matrix: from the standard basis $\mathcal{E}=\{e_1, e_2, e_3, e_4, e_5, e_6 \}$ into the new basis $\mathcal{D}$, see Section~\ref{se:C1} for details.  In this basis, the diagonal matrices, the group of upper triangular matrices and group of unipotent upper triangular matrices in $\Sp_6(\Omega'')$ form a maximal torus, a Borel subgroup $\mathrm{B}$ and the unipotent radical $\rU$ of $\rB$ respectively. 

We find that $\rU$ is a nilpotent subgroup of $\mathrm{GL}_6(\R)$. Now, following~\cite[Thm.~2.1]{R72}, if $\Gamma \cap \rU(\Z)$ is a Zariski dense subgroup of $\rU$ then $\rU/\Gamma \cap \rU(\Z)$ is compact, and hence  $\Gamma \cap \rU(\Z)$ has finite index inside $\rU(\Z).$ Hence, to show that $\Gamma\cap \rU(\Z)$ is of finite index in $\rU(\Z)$, it is enough to show that $\Gamma \cap \rU(\Z)$ is Zariski dense in $\rU$, and for  this it is enough to show that $\Gamma$ contains nontrivial unipotent elements corresponding to each of the positive roots, and the arithmeticity of $\Gamma$ follows from~\cite{Tits}. 

However, following~\cite[Thm.~3.5]{V87}, if $\Gamma$ is a Zariski dense subgroup of $\Sp_6(\Omega)(\Z)$, and intersects the highest and second highest root groups non-trivially, then $\Gamma$ has finite index in $\Sp_6(\Omega)(\Z)$. Hence, to show the arithmeticity of the examples under consideration in this article, following~\cite[Thm.~3.5]{V87}, it is enough to find the unipotent elements corresponding to the highest and second highest roots of $\Sp_6(\Omega'')$. In particular, we show that the Zariski dense subgroup $\Gamma'=X^{-1}\Gamma X$ of $\Sp_6(\Omega'')$ contains some non-trivial elements of the unipotent groups $U_{t_1^{2}}$ and $U_{t_1t_2}$ which establishes the arithmeticity of $\Gamma'$ in $\Sp_{6}(\Omega'')$ following~~\cite[Thm.~3.5]{V87}. Consequently, this proves the arithmeticity of $\Gamma$ in $\Sp_6(\Omega')$.

\subsection{Symplectic form preserved by $\Gamma(f,g)$}\label{se:form} 
We know that $\G(f, g)$ preserves a non-degenerate integral symplectic form $\Omega$ on $\Z^{6}$ and $\G(f, g)\subset \Sp_{\Omega}(\Z)$ is Zariski dense by Theorem~6.5 in~\cite{BH}. For our task, computations of the symplectic form $\Omega$ (which is unique up to scalar) is one of the important steps, and in fact we compute them  explicitly in all the examples discussed in this article. For this purpose, we closely follow the notations and method explained in~\cite{SV,S15S,S17} to prove the arithmeticity of various hypergeometric groups of type $\Sp(4)$. we simply adapt the method to perform the computations for the examples under consideration .  

Let $A$ and $B$ be the companion matrices of $f$ and $g$ respectively. Let $\mathcal{E}=\{e_1, e_2, e_3, e_4, e_5,e_6\}$ be the standard basis vectors of $\Q^{6}$ over $\Q$, and $v$ be the last column vector of $C-I$, where $I$ is the identity matrix. Then,  $Cv =v$.
 Therefore, using the invariance of $\Omega$ under the action of $C$, we get that $v$ is $\Omega$-orthogonal to the vectors $e_1, e_2, e_3, e_4,e_5$ and $\Omega(v, e_6) \neq 0$ (since $\Omega$ is non-degenerate).
 We may now assume that $\Omega(v, e_6) = 1$. It can be easily checked that the set $\mathcal{B}=\{v, Bv, B^2 v, B^3 v, B^4v, B^5 v \}$ (similarly the set $ \{v,Av,A^{2}v,A^{3}v,A^{4}v, A^{5}v\}$) is linearly independent over $\Q$.
 Since $\Omega$ is invariant under the action of $A$, that is, 
 $$\Omega(A^{i} v,A^{j} v) = \Omega(A^{i+1} v,A^{j+1} v),\quad \text{for any} \quad i, j \in  \Z,$$ 
 to determine the symplectic form $\Omega$ on $\Q^{6}$, it is enough to compute $$\Omega(v,A^{j}v),\quad \text{for}\quad j = 0, 1, 2, 3, 4,5 .$$
 Also, since $v$ is $\Omega$-orthogonal to the vectors $e_1, e_2, e_3, e_4, e_5$ and $\Omega(v, e_6) = 1$ (say), we get that $\Omega(v,A^{j}v) $ is the coefficient of $e_6$ in $A^{j}v$.
Since the companion matrix $A$ (resp.~$B$) of $f$ (resp.~$g$) maps $e_i$ to
$e_{i+1}$ for $1 \leq i \leq 5$, to know the symplectic forms $\Omega$ preserved by the
symplectic hypergeometric groups it is enough to find the scalars $\Omega(e_1, e_j)$ for $1 \leq j \leq 6$ , since 
$$\Omega(e_i, e_j) = \Omega(Ae_i,Ae_j) = \Omega(e_{i+1}, e_{j+1})\quad \text{for} \quad 1 \leq i, j \leq 5.$$ 
Combining, all the above information about the symplectic form preserved by $\Gamma(f,g)$, we can now describe the matrix form of $\Omega$, which we denote by the same letter. For simplicity, by writing
\begin{align*}
   & 0=\Omega(e_1, e_1),\quad  b=\Omega(e_1, e_2),\quad c=\Omega(e_1, e_3), \\
   & d=\Omega(e_1, e_4), \quad e=\Omega(e_1, e_5), \quad f=\Omega(e_1, e_6),
\end{align*}
we can write the matrix form of $\Omega$ associated to the group $\G(f,g)$ as follows: 

\resizebox{\textwidth}{!}{
\parbox{\textwidth}{
\begin{equation}\label{eqn:MatrixQ} 
\Omega =  \left(\begin{array}{rrrrrrrrr} 
0 & b & c & d & e & f \\
-b & 0 & b & c & d & e\\
-c & -b & 0 & b & c & d\\
-d & -c & -b & 0 & b & c\\
-e & -d & -c & -b & 0 & b\\
-f & -e & -d & -c & -b & 0\end {array} \right).
\end{equation}
}}


\section{Proof of Theorem~\ref{thm:main}}

\subsection{Arithmeticity of C-1}\label{se:C1} In this case, the parameters are 
$$ \alpha=\left(0,0,0,0,\frac{1}{2},\frac{1}{2}\right) \quad \text{and} \quad \beta=\left(\frac{1}{3},\frac{1}{3},\frac{2}{3},\frac{2}{3},\frac{1}{6},\frac{5}{6}\right).$$ 
The corresponding polynomials are $f(x)={x}^{6}-2\,{x}^{5}-{x}^{4}+4\,{x}^{3}-{x}^{2}-2\,x+1$  and $g(x)={x}^{6}+{x}^{5}+2\,{x}^{4}+{x}^{3}+2\,{x}^{2}+x+1.$
Therefore, $f(x)-g(x)=-3\,{x}^{5}-3\,{x}^{4}+3\,{x}^{3}-3\,{x}^{2}-3\,x$. Let $A$ and $B$ be the companion matrices of $f(x)$ and $g(x)$ respectively. Then
 
\resizebox{\textwidth}{!}{
\parbox{\textwidth}{
\begin{align*}
A=\left(\begin{array}{rrrrrrrr} 
0&0&0&0&0&-1\\ 
1&0&0&0&0&2\\ 
0&1&0&0&0&1\\ 
0&0&1&0&0&-4\\ 
0&0&0&1&0&1\\ 
0&0&0&0&1&2\end {array} \right),\quad 
B=\left(\begin {array}{rrrrrrrr} 
0&0&0&0&0&-1\\ 
1&0&0&0&0&-1\\ 
0&1&0&0&0&-2\\ 
0&0&1&0&0&-1\\ 
0&0&0&1&0&-2\\
0&0&0&0&1&-1
\end {array} \right),\quad 
C=A^{-1}B=\left(\begin{array}{rrrrrrrr} 
1&0&0&0&0&-3\\ 
0&1&0&0&0&-3\\ 
0&0&1&0&0&3\\ 
0&0&0&1&0&-3\\ 
0&0&0&0&1&-3\\ 
0&0&0&0&0&1
\end {array} \right).
\end{align*}}}
The hypergeometric group $\Gamma(f,g)=\langle A, B \rangle$ is a subgroup of $\mathrm{SL}_6(\Z)$, preserves a symplectic form $\Omega$ which we now compute following the discussion of Section~\ref{se:form}. 

Let $\mathcal{E}=\{ e_1, e_2, e_3, e_4, e_5, e_6\}$ be the standard basis of $\Q^{6}$ over $\Q$, and let $v=(C-I)e_6$. Then 
\begin{align*}
    & v= -3e_1-3e_2+3e_3-3e_4-3e_5,\\
    & B v= -3e_2-3e_3+3e_4-3e_5-3e_6,\\
    & B^2 v= 3e_1+3e_2+3e_3+9e_5,\\
    & B^3 v= 3e_2+3e_3+3e_4+9e_6,\\
    & B^4 v= -9e_1-9e_2-15e_3-6e_4-15e_5-9e_6,\\
    & B^5 v=9e_1+9e_3-6e_4+12e_5-6e_6.
\end{align*}
As mentioned, the set $\mathcal{B}=\{v, Bv, B^2 v, B^3 v, B^4v, B^5 v \}$ forms a basis of $\Q^6$, and we find that with respect to the basis $\mathcal{B}$, the matrix of the symplectic form preserved by $\Gamma$ is given by 

\resizebox{\textwidth}{!}{
\parbox{\textwidth}{
\begin{align*}
\Omega=\left(\begin{array}{cccccc} 0&-3&0&9&-9&-6\\ \noalign{\medskip}3&0&
-3&0&9&-9\\ \noalign{\medskip}0&3&0&-3&0&9\\ \noalign{\medskip}-9&0&3&0
&-3&0\\ \noalign{\medskip}9&-9&0&3&0&-3\\ \noalign{\medskip}6&9&-9&0&3
&0\end {array} \right)
\end{align*}
}}

Now, let $Z$ be the change of basis matrix from the basis $\mathcal{B}$ to the standard basis $\mathcal{E}$. Then, following a simple computation we obtain that

\resizebox{\textwidth}{!}{
\parbox{\textwidth}{
\begin{align*}
Z=\left(\begin{array}{rrrrrrrrrr} 
-{\frac {5}{27}}&-{\frac {4}{27}}& \frac{1}{9}&\frac{1}{27}&{\frac {2}{27}}&0\\ \noalign{\medskip}
\frac{1}{27}&-\frac{1}{3}&-\frac{1}{27}&{\frac {4}{27}}&\frac{1}{9}&{\frac {2}{27}}\\ \noalign{\medskip}
-{\frac {5}{27}}&-{\frac{7}{27}}&-\frac{1}{9}&\frac{1}{27}&{\frac {8}{27}}&\frac{1}{9}\\ \noalign{\medskip}
{\frac {4}{27}}&-\frac{1}{3}&-{\frac {4}{27}}&-{\frac {2}{27}}&\frac{1}{9}&{\frac {8}{27}}\\ \noalign{\medskip}
\frac{1}{27}&-{\frac {4}{27}}&-\frac{1}{9}&-{\frac {2}{27}}&{\frac {2}{27}}&\frac{1}{9}\\ \noalign{\medskip}
{\frac {4}{27}}&-\frac{1}{9}&-\frac{1}{27}&-{\frac {2}{27}}&0&{\frac {2}{27}}
\end {array} \right)
\end{align*}
}}
 and with respect to the standard basis $\mathcal{E}$, the matrix, up to scalar multiplication, of the symplectic form preserved by the hypergeometric group $\Gamma$ is given by 
 
\resizebox{\textwidth}{!}{
\parbox{\textwidth}{
\begin{align*}
\Omega'=27 Z^{t}\Omega Z=\left(\begin {array}{rrrrrrrr} 
0&2&1&3&-4&-5\\ \noalign{\medskip}
-2&0&2&1&3&-4\\ \noalign{\medskip}
-1&-2&0&2&1&3\\ \noalign{\medskip}
-3&-1&-2&0&2&1\\ \noalign{\medskip}
4&-3&-1&-2&0&2\\ \noalign{\medskip}
5&4&-3&-1&-2&0\end {array} \right)\,
\end{align*}
}}

where $Z^{t}$ denotes the transpose of the matrix $Z$ and multiplication by 27 is purely aesthetic purposes, that is to have cleaner entries in the matrix of $\Omega$. We now find that $A^t \Omega'  A = \Omega' = B^t \Omega'  B$.

We now change the standard basis $\mathcal{E}$ to $\mathcal{D}$ for which the change of basis matrix takes the following shape

\resizebox{\textwidth}{!}{
\parbox{\textwidth}{
\begin{align*}
X=\left(\begin {array}{rrrrrrrrr} 
0&0&0&0&0&-3\\ \noalign{\medskip}
-3&-6&0&0&0&-3\\ \noalign{\medskip}
-3&-6&0&0&6&3\\ \noalign{\medskip}
3&12&-{\frac{27}{2}}&0&-3&-3\\ \noalign{\medskip}
-3&-9&-{\frac{27}{2}}&{\frac{27}{2}}&6&-3\\ \noalign{\medskip}
-3&0&0&0&0&0\end{array}\right).
\end{align*}
}}

Now, with respect to this new  basis $\mathcal{D}$, we find that $\Omega'$ takes the anti diagonal form denoted by $\Omega''=X^{t}\Omega' X$, that is,
 
\resizebox{\textwidth}{!}{
\parbox{\textwidth}{
\begin{align*}
\Omega''=\left(\begin{array}{rrrrrrrrr}
0&0&0&0&0&81\\ \noalign{\medskip}
0&0&0&0&-162&0\\ \noalign{\medskip}
0&0&0&{-\frac{729}{2}}&0&0\\ \noalign{\medskip}
0&0&{\frac {729}{2}}&0&0&0\\ \noalign{\medskip}
0&162&0&0&0&0\\ \noalign{\medskip}
-81&0&0&0&0&0\end {array} \right).
\end{align*}
}}

With respect to the basis $\mathcal{D}$, we write the generators $a=X^{-1}AX, b=X^{-1}BX$ and $c=a^{-1}b$ as follows:
  
\resizebox{\textwidth}{!}{
\parbox{\textwidth}{
\begin{align*}
a=\left(\begin{array}{rrrrrrrrrr} 
3&3&9/2&-9/2&-2&1\\ \noalign{\medskip}0
&-3/2&-9/4&9/4&1&0\\ \noalign{\medskip}
0&0&-1&1&0&0\\ \noalign{\medskip}
0&1&-5/2&3/2&0&0\\ \noalign{\medskip}
1&-1&0&0&0&0\\ \noalign{\medskip}
-1&0&0&0&0&0\end{array} \right)
\quad 
b=\left(\begin{array}{rrrrrrrrr}
0&3&9/2&-9/2&-2&1\\ \noalign{\medskip}
0&-3/2&-9/4&9/4&1&0\\ \noalign{\medskip}
0&0&-1&1&0&0\\ \noalign{\medskip}
0&1&-5/2&3/2&0&0\\ \noalign{\medskip}
1&-1&0&0&0&0\\ \noalign{\medskip}
-1&0&0&0&0&0\end{array}\right) \,, 
\quad c= \left( \begin {array}{rrrrrrrr} 
1&0&0&0&0&0\\ \noalign{\medskip}
0&1&0&0&0&0\\ \noalign{\medskip}
0&0&1&0&0&0\\ \noalign{\medskip}
0&0&0&1&0&0\\ \noalign{\medskip}
0&0&0&0&1&0\\ \noalign{\medskip}
-3&0&0&0&0&1\end {array} \right)\,.
\end{align*}}}

Now, note that the elements $c=a^{-1}b$   and $ba^{-1}$ in these cases are always  unipotent elements. Therefore, in this particular case, we write $q_1=b a^{-1}$ and consider the elements
\begin{align*}
& w_1=[a \, b^{-1}], \quad w_2=[a\, b^2], \quad w_3=w_2 q_1^{-8}, \quad w_4=q_1 c\\
& w_5=[w_1\, w_3], \quad w_6=[w_4 \, w_3], \quad w_7= [w_6\, w_5].
\end{align*}
Following simple computation, we find that the $q_1$ and $q_2=w_3^{2916} w_{7}^{-1}$ are the desired unipotent elements corresponding to the highest and the second highest roots of $\Sp_6$:

\resizebox{\textwidth}{!}{
\parbox{\textwidth}{
 \begin{align*}
q_1=b a^{-1} =\left(\begin {array}{cccccc} 1&0&0&0&0&3\\ 0&1&0&0
&0&0\\ 0&0&1&0&0&0\\ 0&0&0&1&0&0
\\ 0&0&0&0&1&0\\ 0&0&0&0&0&1
\end {array} \right) \,, \quad q_2= \left( \begin {array}{cccccc} 1&0&0&0&52488&0\\ 0&1
&0&0&0&-26244\\ 0&0&1&0&0&0\\ 0&0&0
&1&0&0\\ 0&0&0&0&1&0\\ 0&0&0&0&0&1
\end {array}  \right).\end{align*}}}

Hence, following the discussion in Section~\ref{se:meth}, the arithmeticity follows by the existence of the unipotents $q_1$ and $q_2$. 

\subsection{Arithmeticity of C-10}\label{se:C10} In this case, the parameters are
$$ \left(0,0,0,0,\frac{1}{3},\frac{2}{3}\right) \quad\text{and}\quad \left(\frac{1}{9},\frac{2}{9},\frac{4}{9},\frac{5}{9},\frac{7}{9},\frac{8}{9}\right).$$ The corresponding polynomials are $f(x)={x}^{6}-3\,{x}^{5}+3\,{x}^{4}-2\,{x}^{3}+3\,{x}^{2}-3\,x+1$ and $g(x)={x}^{6}+{x}^{3}+1.$ Therefore, $f(x)-g(x)=-3\,{x}^{5}+3\,{x}^{4}-3\,{x}^{3}+3\,{x}^{2}-3\,x$ and

\resizebox{\textwidth}{!}{
\parbox{\textwidth}{
\begin{align*}
A=\left(\begin{array}{rrrrrrrrr} 
0&0&0&0&0&-1\\ \noalign{\medskip}
1&0&0&0&0&3\\ \noalign{\medskip}
0&1&0&0&0&-3\\ \noalign{\medskip}
0&0&1&0&0&2\\ \noalign{\medskip}
0&0&0&1&0&-3\\ \noalign{\medskip}
0&0&0&0&1&3\end {array} \right),
\quad 
B=\left(\begin{array}{rrrrrrrr} 
0&0&0&0&0&-1\\ \noalign{\medskip}
1&0&0&0&0&0\\ \noalign{\medskip}
0&1&0&0&0&0\\ \noalign{\medskip}
0&0&1&0&0&-1\\ \noalign{\medskip}
0&0&0&1&0&0\\ \noalign{\medskip}
0&0&0&0&1&0\end {array} \right), 
\quad
C= A^{-1}B=\left(\begin{array}{rrrrrrrr}
1&0&0&0&0&-3\\ \noalign{\medskip}
0&1&0&0&0&3\\ \noalign{\medskip}
0&0&1&0&0&-3\\ \noalign{\medskip}
0&0&0&1&0&3\\ \noalign{\medskip}
0&0&0&0&1&-3\\ \noalign{\medskip}
0&0&0&0&0&1\end{array}\right).
\end{align*}}}

Now, we write down the symplectic form $\Omega'$ (up to scalar multiple) preserved by $\Gamma=\langle A, B\rangle$, with respect to the basis $\mathcal{D}$ the change of basis matrix $X$, and symplectic form $\Omega''$ preserved by the group $\Gamma'=X^{-1} \Gamma X = \langle a=X^{-1}AX, b=X^{-1}BX\rangle :$

\resizebox{\textwidth}{!}{
\parbox{\textwidth}{
\begin{align*}
\Omega'=\left(\begin {array}{rrrrrrrr} 0&1&2&2&1&-1\\ 
-1&0&1&2&2&1\\ \noalign{\medskip}
-2&-1&0&1&2&2\\ \noalign{\medskip}
-2&-2&-1&0&1&2\\ \noalign{\medskip}
-1&-2&-2&-1&0&1\\ \noalign{\medskip}
1&-1&-2&-2&-1&0\end {array} \right)\,,
\quad 
X=\left(\begin {array}{rrrrrrrrr} 
0&0&0&0&0&-1\\ \noalign{\medskip}
-1&1&0&0&0&1\\ \noalign{\medskip}
1&0&0&0&-1&-1\\ \noalign{\medskip}
-1&-2&-1&0&2&1\\ \noalign{\medskip}
1&2&0&1&-1&-1\\ \noalign{\medskip}
-1&0&0&0&0&0\end {array} \right),
\quad 
\Omega''=\left(\begin {array}{rrrrrrrr} 
0&0&0&0&0&1\\ \noalign{\medskip}
0&0&0&0&1&0\\ \noalign{\medskip}
0&0&0&-1&0&0\\ \noalign{\medskip}
0&0&1&0&0&0\\ \noalign{\medskip}
0&-1&0&0&0&0\\ \noalign{\medskip}
-1&0&0&0&0&0\end{array}\right)\,,
\end{align*}}}

where

\resizebox{\textwidth}{!}{
\parbox{\textwidth}{
\begin{align*}
a= \left( \begin {array}{rrrrrrrr} 
2&-2&0&-1&1&1\\ \noalign{\medskip}0&-2&0
&-1&1&0\\ \noalign{\medskip}0&0&0&1&0&0\\ \noalign{\medskip}0&1&-1&2&0
&0\\ \noalign{\medskip}1&-3&0&-1&1&0\\ \noalign{\medskip}-1&0&0&0&0&0
\end {array} \right)\,,\qquad 
b= \left( \begin {array}{rrrrrrrr} 
-1&-2&0&-1&1&1\\ \noalign{\medskip}0&-2
&0&-1&1&0\\ \noalign{\medskip}0&0&0&1&0&0\\ \noalign{\medskip}0&1&-1&2
&0&0\\ \noalign{\medskip}1&-3&0&-1&1&0\\ \noalign{\medskip}-1&0&0&0&0&0
\end {array} \right)\,,
\quad
c=a^{-1}b=\left(\begin{array}{rrrrrrrr} 1&0&0&0&0&0\\ \noalign{\medskip}0&1&0&0
&0&0\\ \noalign{\medskip}0&0&1&0&0&0\\ \noalign{\medskip}0&0&0&1&0&0
\\ \noalign{\medskip}0&0&0&0&1&0\\ \noalign{\medskip}-3&0&0&0&0&1
\end{array}\right).
\end{align*}}}

In this particular case, we write $q_1=b a^{-1}$ which forms a unipotent element with respect to highest root. Then consider the elements
\begin{align*}
&w_1=[a\,\, b^2], \quad w_2=c^{-1}\left(b^{2} a\right)^{-3} q_1, \quad w_3=c w_2\\
& w_4= [w_1 \, \, w_3], \quad w_5=[w_4 \,\,w_2], \quad w_6= w_3^{59} w_5 .
\end{align*}

This gives us the desired unipotent element $q_2=w_6 c$ with respect to the second highest root.  Hence, the arithmeticity follows by the existence of the unipotents $q_1$ and $q_2$.

\resizebox{\textwidth}{!}{
\parbox{\textwidth}{
\begin{align*}
q_1=b a^{-1} =\left(\begin{array}{rrrrrrr} 
1&0&0&0&0&3\\ \noalign{\medskip}
0&1&0&0&0&0\\\noalign{\medskip} 
0&0&1&0&0&0\\\noalign{\medskip}
0&0&0&1&0&0\\\noalign{\medskip}
0&0&0&0&1&0\\\noalign{\medskip}
0&0&0&0&0&1\end {array} \right) \,,
\quad 
q_2= \left(\begin{array}{rrrrrrrr}
1&0&0&0&-198&0\\\noalign{\medskip} 
0&1&0&0&0&-198\\\noalign{\medskip}
0&0&1&0&0&0\\\noalign{\medskip}
0&0&0&1&0&0\\\noalign{\medskip}
0&0&0&0&1&0\\\noalign{\medskip}
0&0&0&0&0&1\end {array}\right).\end{align*}}}

\subsection{Arithmeticity of C-42} In this case,  the parameters are $$\left(0,0,\frac{1}{4},\frac{1}{4},\frac{3}{4},\frac{3}{4}\right) \quad \text{and} \quad  \left(\frac{1}{3},\frac{2}{3},\frac{1}{12},\frac{5}{12},\frac{7}{12},\frac{11}{12}\right).$$
The corresponding polynomials are 
$f(x)={x}^{6}-2\,{x}^{5}+3\,{x}^{4}-4\,{x}^{3}+3\,{x}^{2}-2\,x+1$ and $g(x)={x}^{6}+{x}^{5}-{x}^{3}+x+1.$ Therefore, $f(x)-g(x)=-3\,{x}^{5}+3\,{x}^{4}-3\,{x}^{3}+3\,{x}^{2}-3\,x$ and

\resizebox{\textwidth}{!}{
\parbox{\textwidth}{
\begin{align*}
A= \left( \begin {array}{rrrrrrrr} 0&0&0&0&0&-1\\ \noalign{\medskip}1&0&0&0
&0&2\\ \noalign{\medskip}0&1&0&0&0&-3\\ \noalign{\medskip}0&0&1&0&0&4
\\ \noalign{\medskip}0&0&0&1&0&-3\\ \noalign{\medskip}0&0&0&0&1&2
\end {array} \right), 
\quad B= \left( \begin {array}{rrrrrrrrr} 0&0&0&0&0&-1\\ \noalign{\medskip}1&0&0&0
&0&-1\\ \noalign{\medskip}0&1&0&0&0&0\\ \noalign{\medskip}0&0&1&0&0&1
\\ \noalign{\medskip}0&0&0&1&0&0\\ \noalign{\medskip}0&0&0&0&1&-1
\end {array} \right),
\quad 
C=A^{-1}B=\left(\begin{array}{rrrrrrrrr} 1&0&0&0&0&-3\\ \noalign{\medskip}0&1&0&0
&0&3\\ \noalign{\medskip}0&0&1&0&0&-3\\ \noalign{\medskip}0&0&0&1&0&3
\\ \noalign{\medskip}0&0&0&0&1&-3\\ \noalign{\medskip}0&0&0&0&0&1
\end {array} \right).
\end{align*}}}
Now, we write down the symplectic form $\Omega'$ (up to scalar multiple) preserved by $\Gamma=\langle A, B\rangle$, with respect to the basis $\mathcal{D}$ the change of basis matrix $X$, and symplectic form $\Omega''$ preserved by the group $\Gamma'=X^{-1} \Gamma X = \langle a=X^{-1}AX, b=X^{-1}BX\rangle :$

\resizebox{\textwidth}{!}{
\parbox{\textwidth}{
\begin{align*}
\Omega'=\left( \begin {array}{cccccc} 0&0&-1&-1&0&-1\\ \noalign{\medskip}0&0&0
&-1&-1&0\\ \noalign{\medskip}1&0&0&0&-1&-1\\ \noalign{\medskip}1&1&0&0
&0&-1\\ \noalign{\medskip}0&1&1&0&0&0\\ \noalign{\medskip}1&0&1&1&0&0
\end {array} \right)\,, 
\quad X=\left( \begin {array}{cccccc} 0&0&0&-3&0&0\\ \noalign{\medskip}0&3&-3
&3&0&3/2\\ \noalign{\medskip}0&-3&3&-3&0&3/2\\ \noalign{\medskip}0&-6&
-3&3&6&0\\ \noalign{\medskip}-3/2&3&3&-3&-3&-3/2\\ \noalign{\medskip}0
&0&-3&0&0&0\end {array} \right)\,,
\quad
\Omega''= \left( \begin {array}{cccccc} 0&0&0&0&0&-9/2\\ \noalign{\medskip}0&0&0
&0&-18&0\\ \noalign{\medskip}0&0&0&9&0&0\\ \noalign{\medskip}0&0&-9&0&0
&0\\ \noalign{\medskip}0&18&0&0&0&0\\ \noalign{\medskip}9/2&0&0&0&0&0
\end {array} \right),
\end{align*}
}}
where

\resizebox{\textwidth}{!}{
\parbox{\textwidth}{
\begin{align*}
a = \left( \begin {array}{cccccc} 1/2&3&1&0&-3&-1/2\\ \noalign{\medskip}1
/2&-3/2&0&0&1&1/4\\ \noalign{\medskip}1/2&-1&1&1&1&1/2
\\ \noalign{\medskip}0&0&-1&0&0&0\\ \noalign{\medskip}3/4&-5/2&-1/2&0&
3/2&3/4\\ \noalign{\medskip}0&1&0&0&0&1/2\end {array} \right)
\,, \quad b=\left( \begin {array}{cccccc} 1/2&3&1&0&-3&-1/2\\ \noalign{\medskip}1
/2&-3/2&0&0&1&1/4\\ \noalign{\medskip}1/2&-1&-2&1&1&1/2
\\ \noalign{\medskip}0&0&-1&0&0&0\\ \noalign{\medskip}3/4&-5/2&-1/2&0&
3/2&3/4\\ \noalign{\medskip}0&1&0&0&0&1/2\end {array} \right)\,,
\quad 
c=a^{-1}b= \left( \begin {array}{cccccc} 1&0&0&0&0&0\\ \noalign{\medskip}0&1&0&0
&0&0\\ \noalign{\medskip}0&0&1&0&0&0\\ \noalign{\medskip}0&0&-3&1&0&0
\\ \noalign{\medskip}0&0&0&0&1&0\\ \noalign{\medskip}0&0&0&0&0&1
\end {array} \right).
\end{align*}
}}

Then consider the elements
\begin{align*} \begin {array}{lllllll} 
    & w_1=[a\, b^{-1}] & w_2= a^3 c a^{-3} & w_3=a^{4} c a^{-4} & w_4 = b^3 c b^{-3} & w_5=w_2 c^{-1} \\\noalign{\medskip}
    & w_6=[w_1 \, w_5] & w_7=w_{5}^{4} & w_8 = b^2 c b^{-2} & w_9=[w_4\, w_5] &  w_{10}=w_{9}^{4} w_{7}^{-18} \\\noalign{\medskip}
    & w_{11}= w_{6}^{4} & w_{12} = w_{11}^{9} w_{10} & w_{13}=[w_5\, w_8] & w_{14}=w_{13}^{2} w_{7}^{3} w_{11}^{-2} & w_{15}= w_{14}^{1938} w_{12}^{-73} \\\noalign{\medskip}
    & w_{16} =[w_3\, w_{15}] & w_{17}=w_{15}^{9} w_{16} & w_{18}= w_3 w_{15} w_{3}^{-1} w_{17}^{-1} .
    \end {array}
\end{align*}
This will give us the desired unipotent elements $$q_1=(w_{15}^{8} w_{18})^{-1}\qquad \text{and} \qquad q_2= w_{17}^{-15570} q_1^{77849}, $$ namely, 

\resizebox{\textwidth}{!}{
\parbox{\textwidth}{
\begin{align*}
q_1=\left( \begin {array}{cccccc} 1&0&0&0&0&104727556800
\\ \noalign{\medskip}0&1&0&0&0&0\\ \noalign{\medskip}0&0&1&0&0&0
\\ \noalign{\medskip}0&0&0&1&0&0\\ \noalign{\medskip}0&0&0&0&1&0
\\ \noalign{\medskip}0&0&0&0&0&1\end {array} \right) \,, \qquad q_2= \left( \begin {array}{cccccc} 1&0&0&0&17454592800&0
\\ \noalign{\medskip}0&1&0&0&0&4363648200\\ \noalign{\medskip}0&0&1&0&0
&0\\ \noalign{\medskip}0&0&0&1&0&0\\ \noalign{\medskip}0&0&0&0&1&0
\\ \noalign{\medskip}0&0&0&0&0&1\end {array} \right)\,. 
\end{align*}
}}

Hence, the arithmeticity follows by the existence of the unipotents $q_1$ and $q_2$.

\subsection{Arithmeticity of C-59} In this case, the parameters are
$$ \left(0,0,\frac{1}{12},\frac{5}{12},\frac{7}{12},\frac{11}{12}\right) \quad \text{and} \quad \left(\frac{1}{3},\frac{2}{3},\frac{1}{4},\frac{3}{4},\frac{1}{4},\frac{3}{4}\right) .$$ 

The corresponding polynomials are $f(x)={x}^{6}-2\,{x}^{5}+2\,{x}^{3}-2\,x+1$ and  $g(x)={x}^{6}+{x}^{5}+3\,{x}^{4}+2\,{x}^{3}+3\,{x}^{2}+x+1.$ Therefore, $f(x)-g(x)= -3 x^5 -3 x^4-3x^2-3x$, and

\resizebox{\textwidth}{!}{
\parbox{\textwidth}{
\begin{align*}
A= \left( \begin {array}{rrrrrrrr} 
0&0&0&0&0&-1\\ \noalign{\medskip}1&0&0&0
&0&2\\ \noalign{\medskip}0&1&0&0&0&0\\ \noalign{\medskip}0&0&1&0&0&-2
\\ \noalign{\medskip}0&0&0&1&0&0\\ \noalign{\medskip}0&0&0&0&1&2
\end {array} \right), 
\quad B= \left( \begin {array}{rrrrrrrrr} 
0&0&0&0&0&-1\\ \noalign{\medskip}1&0&0&0
&0&-1\\ \noalign{\medskip}0&1&0&0&0&-3\\ \noalign{\medskip}0&0&1&0&0&-
2\\ \noalign{\medskip}0&0&0&1&0&-3\\ \noalign{\medskip}0&0&0&0&1&-1
\end {array} \right),
\quad 
C=A^{-1}B=\left(\begin{array}{rrrrrrrrr} 
1&0&0&0&0&-3\\ \noalign{\medskip}0&1&0&0
&0&-3\\ \noalign{\medskip}0&0&1&0&0&0\\ \noalign{\medskip}0&0&0&1&0&-3
\\ \noalign{\medskip}0&0&0&0&1&-3\\ \noalign{\medskip}0&0&0&0&0&1
\end {array} \right).
\end{align*}}}
Now, we write down the symplectic form $\Omega'$ (up to scalar multiple) preserved by $\Gamma=\langle A, B\rangle$, with respect to the basis $\mathcal{D}$ the change of basis matrix $X$ and symplectic form $\Omega''$ preserved by the group $\Gamma'=X^{-1} \Gamma X = \langle a=X^{-1}AX, b=X^{-1}BX\rangle :$

\resizebox{\textwidth}{!}{
\parbox{\textwidth}{
\begin{align*}
\Omega'=\left( \begin {array}{cccccc} 0&1&0&1&-2&-3\\ \noalign{\medskip}-1&0&
1&0&1&-2\\ \noalign{\medskip}0&-1&0&1&0&1\\ \noalign{\medskip}-1&0&-1&0
&1&0\\ \noalign{\medskip}2&-1&0&-1&0&1\\ \noalign{\medskip}3&2&-1&0&-1
&0\end {array} \right)\,, 
\quad 
X=\left( \begin {array}{cccccc} -3&0&0&0&0&0\\ \noalign{\medskip}-3&0&0
&0&-6&-3\\ \noalign{\medskip}0&6&0&0&12&-3\\ \noalign{\medskip}-3&0&0&
-12&6&0\\ \noalign{\medskip}-3&6&12&24&12&-3\\ \noalign{\medskip}0&0&0
&0&0&-3\end {array} \right)\,,
\quad
\Omega''=\left( \begin {array}{cccccc} 0&0&0&0&0&-36\\ \noalign{\medskip}0&0&0
&0&72&0\\ \noalign{\medskip}0&0&0&144&0&0\\ \noalign{\medskip}0&0&-144
&0&0&0\\ \noalign{\medskip}0&-72&0&0&0&0\\ \noalign{\medskip}36&0&0&0&0
&0\end {array} \right)\,,
\end{align*}}}
where

\resizebox{\textwidth}{!}{
\parbox{\textwidth}{
\begin{align*}
a = \left( \begin {array}{rrrrrrrr} 0&0&0&0&0&-1\\ \noalign{\medskip}0&-3&-
6&-12&-7&1\\ \noalign{\medskip}0&0&-2&-5&1&0\\ \noalign{\medskip}0&0&1
&2&0&0\\ \noalign{\medskip}0&1&2&4&2&0\\ \noalign{\medskip}1&-2&-4&-8&
-4&3\end {array} 
\right)\,, \qquad b = \left(\begin {array}{rrrrrrrr} 0&0&0&0&0&-1\\ \noalign{\medskip}0&-3&-
6&-12&-7&1\\ \noalign{\medskip}0&0&-2&-5&1&0\\ \noalign{\medskip}0&0&1
&2&0&0\\ \noalign{\medskip}0&1&2&4&2&0\\ \noalign{\medskip}1&-2&-4&-8&
-4&0\end {array} \right)\,,
\quad
c=a^{-1}b=\left(\begin {array}{rrrrrrr} 1&0&0&0&0&-3\\ \noalign{\medskip}0&1&0&0
&0&0\\ \noalign{\medskip}0&0&1&0&0&0\\ \noalign{\medskip}0&0&0&1&0&0
\\ \noalign{\medskip}0&0&0&0&1&0\\ \noalign{\medskip}0&0&0&0&0&1
\end {array} \right). 
\end{align*}}}

Then consider the elements
\begin{align*}
    & w_1= [a \,\, b]\,, \quad w_2= [a \,\, b^{-1}]\,, \quad w_3=[b^2 \,\, a^{-1}]\,,  \quad w_4 =[w_1 \,\, w_2] \,, \quad w_5= [w_3 \,\, w_4]  \,.
\end{align*}

This will give us the desired unipotent elements $q_1=c$ and $q_2=w_3^{180} w_5$.

\resizebox{\textwidth}{!}{
\parbox{\textwidth}{
\begin{align*}
q_1=c=\left(\begin{array}{rrrrrrrr} 
1&0&0&0&0&-3\\ \noalign{\medskip}
0&1&0&0&0&0\\ \noalign{\medskip}
0&0&1&0&0&0\\ \noalign{\medskip}
0&0&0&1&0&0\\ \noalign{\medskip}
0&0&0&0&1&0\\ \noalign{\medskip}
0&0&0&0&0&1\end {array}\right), 
\qquad 
q_2=\left(\begin{array}{rrrrrrrr} 
1&0&0&0&-1080&0\\ \noalign{\medskip}
0&1&0&0&0&540\\ \noalign{\medskip}
0&0&1&0&0&0\\ \noalign{\medskip}
0&0&0&1&0&0\\ \noalign{\medskip}
0&0&0&0&1&0\\ \noalign{\medskip}
0&0&0&0&0&1\end{array}\right).
\end{align*}}}

Hence, the arithmeticity follows by the existence of the unipotents $q_1$ and $q_2$.

\section*{Acknowledgements} 
The author would like to thank the Max Planck Institute f\"ur Mathematics (MPIM), Bonn where much of the work on this article was accomplished, for the hospitality and support.


\nocite{}
\bibliographystyle{abbrv}
\bibliography{B}
\end{document}